\documentclass[12pt]{amsart}
\usepackage{amscd,amsmath,amssymb,amsfonts}
\usepackage[cmtip, all]{xy}
\theoremstyle{plain}
\newtheorem{thm}{Theorem}
\newtheorem{lem}[thm]{Lemma}

\newtheorem{prop}[thm]{Proposition}

\theoremstyle{definition}

\newtheorem{ex}[thm]{Example}

\numberwithin{thm}{section}
\numberwithin{equation}{section}

\newcommand{\eq}[2]{\begin{equation}\label{#1}#2 \end{equation}}
\newcommand{\ml}[2]{\begin{multline}\label{#1}#2 \end{multline}}
\newcommand{\ga}[2]{\begin{gather}\label{#1}#2 \end{gather}}

\newcommand{\surj}{\twoheadrightarrow}
\newcommand{\inj}{\hookrightarrow}

\newcommand{\codim}{{\rm codim}}

\newcommand{\Spec}{{\rm Spec \,}}

\newcommand{\sD}{{\mathcal D}}

\newcommand{\sI}{{\mathcal I}}

\newcommand{\sO}{{\mathcal O}}
\newcommand{\sP}{{\mathcal P}}
\newcommand{\sQ}{{\mathcal Q}}

\newcommand{\sU}{{\mathcal U}}
\newcommand{\sV}{{\mathcal V}}

\newcommand{\sZ}{{\mathcal Z}}
\newcommand{\A}{{\mathbb A}}

\newcommand{\C}{{\mathbb C}}

\newcommand{\F}{{\mathbb F}}

\newcommand{\N}{{\mathbb N}}
\renewcommand{\P}{{\mathbb P}}
\newcommand{\Q}{{\mathbb Q}}

\newcommand{\Z}{{\mathbb Z}}

\newcommand{\0}{\emptyset}
\newcommand{\id}{{\rm id}}
\newcommand{\cl}{{\rm cl}}
\begin{document}

\title[Decomposition]{Decomposition of the diagonal and eigenvalues of
Frobenius for Fano hypersurfaces.
}
\author{Spencer Bloch}
\address{University of Chicago, Mathematics, IL 60 636, Chicago, USA}
\email{bloch@math.uchicago.edu}
\author{H\'el\`ene Esnault}
\address{Mathematik,
Universit\"at Essen, FB6, Mathematik, 45117 Essen, Germany}
\email{esnault@uni-essen.de}
\author{Marc Levine}
\address{Northeastern University, Mathematics, MA 02115, Boston, USA  }
\email{marc@neu.edu}
\date{Feb. 10, 2003}
\begin{abstract}
Let $X\subset \P^n$ be a
possibly singular
hypersurface of degree $d\le n$, defined
over a finite field $\F_q$. We show that
the diagonal, suitably interpreted, is decomposable. This gives a proof
that
the eigenvalues
  of the Frobenius action on its $\ell$-adic cohomology
  $H^i(\bar{X}, \Q_\ell)$, for $\ell \neq
{\rm char}(\F_q)$, are divisible by $q$, without using the result
on the existence of rational points by Ax and Katz \cite{Ka}.
\end{abstract}
%\subjclass{Primary Algebraic Geometry}
\maketitle
\begin{quote}

\end{quote}

\section{Introduction}
If $X$ is a  variety defined over a finite field
$k=\F_q$, one encodes the number of its rational points over
all finite extensions $\F_{q^s}\supset \F_q$ in the zeta function,
defined by its logarithmic derivative
\ga{1}{\frac{\zeta'(X, t)}{\zeta(X,t)}=\sum_{s\ge 1} |X(\F_{q^s})|t^{s-1}.}
By the theorem of Dwork \cite{Dw}, we know that
$\zeta(t)$  is a rational function
\ga{2}{\zeta(X,t)\in \Q(t).}
We assume that $X$ is projective and we denote by $U=\P^n\setminus X$ the
complement of a projective embedding.
The Grothendieck-Lefschetz trace formula \cite{Gr} gives a cohomological
formula for the numerator and the denominator of the rational function
\ga{3}{\zeta(U,t)=\prod_{i=0}^{2\ \rm dim(U)} {\rm 
det}(1-F_it)^{(-1)^{i+1}},}
where $F_i$ is the arithmetic Frobenius acting on the compactly supported
$\ell$-adic cohomology  $H^i_c(\overline{U}, \Q_\ell)$.
Letting $H^i_{{\rm prim}}(\overline{X},\Q_\ell)$ denote
the primitive cohomology
$H^i(\overline{X}, \Q_\ell)/H^i(\overline\P^n,\Q_\ell)$ of $X$, we have
\ga{3b}{
H^i_c(\overline{U}, \Q_\ell)
\cong\begin{cases}H^{i-1}_{{\rm prim}}(\overline{X},\Q_\ell)
&\text{ for }(i-1)\le 2 \ {\rm dim}(X), \\
H^i(\overline{ \P^n}, \Q_\ell)&\text{ for }i\ge 2 \ {\rm dim}(X).
\end{cases}
}
For $X$ smooth and complete, the Weil conjectures \cite{DeWeI}
assert that the
  the eigenvalues
of $F_i$ in any complex
embedding $\Q_\ell \subset \C$ have absolute values $q^{\frac{(i-1)}{2}}$.
In particular, there is no possible cancellation of eigenvalues
between the numerator and the denominator of the zeta function.
Consequently,  the property
\ga{4}{|X(\F_{q^s})|\equiv |\P^n(\F_{q^s})| \ \text{mod} \ q^\kappa}
for all $s\ge 1$, and some $\kappa \in \N\setminus \{0\}$,
is equivalent to the property that
\ga{5}{\text{the \ eigenvalues \ of \ } F_i
\text{ are divisible \ by \ }
q^\kappa, \ \\
\text{as \ algebraic \ integers.}\notag}

However, if $X$ is singular, one does not have in general the
  purity of weights of
Frobenius on $H^i_c$. Thus, a   cancellation between the numerator
  of the zeta function and its
denominator is  at least in principle possible, and the property
  \eqref{4} is no longer {\it a
priori} equivalent to the property \eqref{5}.  The purpose of this
   article is to study the
relation between \eqref{4} and \eqref{5} in the   case of
hypersurfaces of degree $d\le n$.

Let $X$ be a complete intersection in $\P^n$
defined by $r$ equations of degrees
$d_1\ge d_2\ge \ldots \ge d_r$ with the property
\ga{6}{ 1\le \kappa=[\frac{n-d_2-\ldots-d_r}{d_1}].}
The theorem of Ax and Katz says precisely that \eqref{4} holds true.
On the other hand, we also know by \cite{DeSGA}, \cite{DD}, \cite{E},
  that if the finite field is replaced by a field of
characteristic 0, the Hodge type of $X$ is $ \kappa$ for all
cohomology groups of $X$.
(See \cite{BE} for a more precise discussion of those theorems).
This gives a strong indication that \eqref{5} should be true as well.
Indeed, as explained to us by Daqing Wan, \eqref{5} is true for $\kappa=1$.
One knows by \cite{DeInt}, Theorem 5.5.3, that $q$ divides the eigenvalues
of Frobenius acting on $H^a(\overline{X}, \Q_\ell)$ for $a> {\rm dim }(X)$,
and since this cohomology vanishes for complete intersections for
$a<{\rm dim}(X)$, the theorem of Ax and Katz implies divisibility
by $q$ for $a={\rm dim}(X)$ as well.
Similarly, using vanishing and Ax-Katz's result, and
replacing \cite{DeInt}, Theorem 5.5.3 by the corresponding
statement for the slopes of the Frobenius action on rigid cohomology
(\cite{W}, p. 820), one obtains that the slopes of the Frobenius action
on rigid cohomology are $\ge 1$.

The purpose of this note is to give a motivic interpretation
for Fano hypersurfaces and $\kappa=1$
of the divisibility result, which does not use the theorem by Ax and Katz.

We now describe our method. Let us first assume that $X$ is smooth.
   By Roitman's theorem
\cite{Ro},  we know that
$CH_0(X\times_k K)=\Z$ for any field extension $K\supset k$ which is
algebraically closed. By \cite[appendix  to lecture 1]{B},
  this implies that the class of the diagonal in $CH_{n-1}(X\times X)$
goes to  zero  in $CH_{n-1}((X\times X\setminus
(\xi\times X \cup X\times A))_\Q$
for some divisor $A$ on $X$ and some 0-cycle $\xi$.
  Letting this class act as a correspondence, it follows that the
  restriction map $H^i(\overline{X}, \Q_\ell)\to
  H^i(\overline{X\setminus A}, \Q_\ell)$
  is zero for $i\ge 1$.
This shows divisibility, as in  \cite{E2}, lemma 2.1.

For singular varieties, the proof of the Hodge-type
statement in the complete intersection
singular case (\cite{E}) shows that the
cohomology with compact support $H_c^i(U)=:H^i(\P^n, X)$
carries the necessary information, and is
easier to deal with than its dual $H^j(U)$.   To carry
out the argument used in the smooth case,
one needs a version of the Chow groups which is related
to compactly supported cohomology. If $X$ is
a strict normal crossing divisor, one can use the
{\em relative motivic cohomology}
$H^{2n}_M(\P^n \times U, X\times U,\Z(n))$,  as defined in
\cite{Le}, chapter 4, 2.2 and p. 209;
this relative motivic cohomology acts as correspondences on
$H^*_c(\overline{U},\Q_\ell)$. Due to the lack of
resolution of singularities in positive
characteristic, we will in general need an {\em alteration}
$\pi: (\P, Y)\to (\P^n, X)$ of
$(\P^n,X)$, that is, a projective, generically finite morphism
$\pi:\P\to\P^n$, with $\P$ smooth,
such that $Y:=\pi^{-1}(X)$ is a strict normal crossing divisor.
We then use the  relative motivic cohomology
$H^{2n}_M(\P \times U, Y\times U,\Z(n))$.

Recall that $H^m_M(\P \times U, Y\times U,\Z(n))$ is the homology
$H_{2n-m}(\sZ^n(\P \times U, Y\times U,*))$, where
  $\sZ^n(\P \times U, Y\times U,*)$ is the single
complex associated to the double
  higher Chow cycle complex
\eq{7}{\begin{CD}
\cdots &  &\cdots & &\cdots\\
@V\partial VV @V\partial VV @VVV\\
\sZ^n(\P\times U, 1)@>{\rm rest}>> \sZ^n(Y^{(1)} \times U, 1)@>{\rm rest}>>
\sZ^n(Y^{(2)}\times U,1)\\
@V\partial VV @V\partial VV @VVV\\
\sZ^n(\P\times U, 0)@>{\rm rest}>> \sZ^N(Y^{(1)} \times U, 0)@>{\rm rest}>>
\sZ^n(Y^{(2)}\times U,0).
\end{CD}}
Here $Y^{(a)}$ is the normalization of all the strata of codimension $a$,
  $\sZ^n(Y^{(a)}\times U, b)$ is a group of cycles on
$Y^{(a)}\times U\times S^b$ where
$S^\bullet$ is the cosimplicial scheme $S^n = \Spec(k[t_0,\dotsc,t_n]/
(\sum t_i -1))$ with
face maps $S^n \inj S^{n+1}$ defined by $t_i=0$.
More precisely, $\sZ^n(Y^{(a)}\times U,
b)$ is generated by the codimension $n$ subvarieties
$Z\subset Y^{(a)}\times U \times S^b$ such
that, for each face $F$ of
$S^b$, and each irreducible component $F'\subset Y^{(a)}$
of the strata of $Y$ we have
$
\codim_{F'\times U\times F}(Z\cap (F'\times U\times F))\ge n.
$
The horizontal restriction maps
are the intersection with the smaller strata, the vertical
$\partial$'s are the boundary maps.

For technical reasons, we find it convenient to use
a subcomplex $\sZ^n(\P \times U, \sI(Y\times U),*)$ of
$\sZ^n(\P \times U, Y\times U,*)$. For $T$ a smooth $k$-scheme
  of finite type,
  and $A$ a closed subset, let $\sZ^n(T,\sI(A),m)$ be the subgroup of
$\sZ^n(T,m)$ consisting of the cycles $W\in \sZ^n(T,m)$
with $\text{Supp}(W)\cap (A\times S^m)=\emptyset$.
The $\sZ^n(T,\sI(A),m)$ evidently form a subcomplex
$\sZ^n(T,\sI(A),*)$  of $\sZ^n(T,*)$, functorial for
  flat pull-back and proper push-forward. Set
\ga{81}{
H^m_M(T, \sI(A),\Z(n)):=
H_{2n-m}(\sZ^n(T, \sI(A),*)).
}
The inclusion
  $\sZ^n(\P \times U, \sI(Y\times U),*)\to \sZ^n(\P \times U,*)$
  extends to a map of complexes
  $\sZ^n(\P \times U, \sI(Y\times U),*)\to \sZ^n(\P \times U, Y\times U,*)$.
  We call $H^m_M(T, \sI(A),\Z(n))$ the {\em motivic cohomology with 
modulus}.

Let $\Delta\subset \P\times U$ be the inverse
image of the diagonal $\subset \P^n\times U$,
i.e., $\Delta$ is the graph of the alteration
$\pi$ restricted to $\P\times U$.
Since ${\rm rest}(\Delta)=0$,$\Delta$ yields a class
\ga{8}{[\Delta]\in H^{2n}_M(\P\times U, \sI(Y\times U),\Z(n)).}
We show
\begin{thm} \label{thm:mot}
Let $X\subset \P^n$ be a hypersurface of degree
$d\le n$ over a  field $k$. Then  there is an alteration
  $\pi: (\P,Y)\to (\P^n,X)$  and a divisor
$A\subset
\P$ which cuts all the strata of $Y$ in codimension $\ge 1$ and such that
the image of $[\Delta]$ in
  $H^{2n}_M((\P\setminus A)\times U, \sI((Y\setminus A)\times U),\Q(n))$
is zero.
\end{thm}
The main idea behind this geometric statement
relies on the following. By a counting argument,
Roitman \cite{Ro} shows that, for a hypersurface
$X\subset \P^n$ of degree $d\le n$,
the correspondence
\ga{9}{\{(x, \ell)\in X\times {\rm Grass}(1,n)\
|\ \ell\subset X\text{ or } \ell \cap X=\{x\}
\text{ for some }x\in X\}
}
dominates $X$.
It follows that the map
$\Z\cong CH_1(\P^n)\to CH_0(X)$ has cokernel
killed by multiplication by $d=\deg X$, where
$CH_0(X)$ is Fulton's homological Chow group
(\cite{Fu}). This implies $CH_0(X)\otimes \Q=\Q.$

We replace Roitman's correspondence
by
\ga{10}{P:=\{(y,\ell)\in \P\times {\rm Grass}(1,n)\ | \
\pi(y)\in \ell, \text{and either }\\\ell\subset X\text{ or }
\ell\cap X=\{x\} \ \text{for \ some} \ x\in X\}.\notag}
We show that $P$ dominates $\P$, and then use the technique
of blowing up strata of $Y$ introduced
in \cite{BJAG} to find the rational equivalence relation
which holds on the complement of some good divisor $A$.
Finally, we show
\begin{thm}\label{thm:coh} Let
$X\subset \P^n$ be a projective variety over a field $k$,
and let $U=\P^n\setminus X$.
Suppose there is an alteration $\pi: (\P, Y)\to (\P^n, X)$
and a divisor $A\subset \P$ which cuts all strata of $Y$in codimension
$\ge 1$, such that
the image of $[\Delta]$ in
  $H^{2n}_M((\P\setminus A)\times U,\sI((Y\setminus A)\times U),\Q(n))$
is zero.
\begin{enumerate}
\item If the characteristic of
the ground field $k$ is 0, then  ${\rm gr}_0^F H^i(X)=0$ for all $i\ge 1$.
\item If $k=\F_q$ is a finite field, then  the eigenvalues of  the
arithmetic Frobenius $F_i$ acting on the compactly supported
$\ell$-adic cohomology  $H^i_c(\overline{U}, \Q_\ell)$
are all divisible by $q$ as algebraic integers for all $i\ge 1$.
\item If $k$ is a perfect field of characteristic $p$, then the slopes
of the Frobenius operator acting on the rigid cohomology
$H^i_{c}(U/K)$ are all $\ge 1$ for all $i\ge 1$.
\end{enumerate}
\end{thm}
To conclude, we remark that this article solves the natural
question posed in the
introduction of \cite{EW},
but only in the case $\kappa=1$. Thinking of
the discussion developed in \cite{BE}, (5.2)  for
  $\kappa\ge 2$ in the smooth case, it is not entirely
clear what the substitute would be for
\ref{thm:mot}. One may also try to generalize
these results for $\kappa = 1$, replacing the
hypersurface $X$ with a more general singular Fano variety.
A singular Fano variety $X$
over a field is a gometrically connected, projective, Cohen-Macaulay
variety such that the reflexive hull of
$\omega_{X_{{\rm reg}}}^N$ is invertible and ample for
some $N \in \Z\setminus \N$. Examples are hypersurfaces
of degree $d\le n$. The question is then whether a Fano variety fulfills
Theorem
\ref{thm:mot}. If yes, as in \cite{E2}, this would
show that a singular Fano variety over a finite field has a rational point.

\noindent {\it Acknowledgements}. We thank
Pierre Berthelot, Pierre Deligne,  V. Srinivas and  in particular
Daqing Wan
for interesting  discussions on topics related to this work.

\section{The proof of Theorem \ref{thm:mot}}
  This section is devoted to the proof of Theorem \ref{thm:mot}.
We fix a base-field $k$ and write $\P^n$ for $\P^n_k$. We want
to show that a certain
class $[\Delta]$ in motivic cohomology with modulus is trivial.
Suppose for a moment we know this
vanishing for $k$ an infinite field. If $k$ is a finite field,
there exist Galois
extensions $k_\ell/k$ with Galois group
$\Z_\ell$ for any prime $\ell$. In particular,
$k_\ell$ is infinite, so $[\Delta_{k_\ell}]=0$ by hypothesis.
  Since the motivic cohomology with modulus
over $k_\ell$ is a direct limit over motivic cohomology with modulus
  over finite subfields $k'
\subset k_\ell$, and since  motivic cohomology with modulus  admits a norm,
  we conclude $[\Delta]$ is killed by
some power of $\ell$.  Since this is true
for  two different $\ell$, and the union of Galois translates
of $A$ in good position with respect to $Y$ is still in good position,
the theorem follows. Thus, we may assume
$k$ is infinite.
In particular, we will use without comment
various general position arguments.

Fix $X\subset
\P^n$ a hypersurface of degree
$d\le n$. We want to define a closed subvariety $Z$ inside
the Grassmann of lines ${\rm Grass}(1,n)$ consisting of lines
``maximally tangent'' to $X$. We have
the incidence correspondence
$\sU := \{(z,\ell)\ |\ z\in \ell \in  {\rm Grass}(1,n)\}$. Define
\eq{2.1}{\sV := \sU \times_{{\rm Grass}(1,n)} \sU =
\{(y,\ell,z)\ |\ y, z\in \ell\}.
}
The $\P^1$-bundle ${\rm pr}_2:\sV \to \sU, (y, \ell, z)\mapsto ( z,\ell)$
has a section $(x,\ell) \mapsto (x,\ell,x)$. Locally on
$\sU$ we may identify $\sV \cong \P^1 \times \sU$ with
homogeneous coordinates $s,t$ in such
a way that the section is given by $t=0$. The section $\sO_{\P^n}
\xrightarrow{X} \sO_{\P^n}(d)$ pulls back to a section
of $p^*\sO(d)_{\P^n}$ under the projection $p: \sV \to \P^n,
(y,\ell,z)\mapsto y$, and the section $X$
restricts to an equation $F(s,t) = F_0s^d+F_1s^{d-t}t+\ldots + F_dt^d$,
where the $F_i$ are
(local) functions on $\sU$. Note $F_0$ is a local defining equation of
$\sU\times_{\P^n} X \subset
\sU$. We are interested in the closed sets defined locally by $F_0=\ldots =
F_{d-1} = 0 \
(\text{ resp. } F_0=\ldots =F_{d} = 0)$.
Denote by
\eq{2.2}{Z_X' \subset Z' \subset {\rm Grass}(1,n)
}
the projection on the Grassmannian of these sets.
Intuitively, an open set $Z^0:= Z'\setminus Z_X'$
consists of lines $d$-fold tangent to $X$ at a point.
Define $Z \subset Z'$ to be the closure of
$Z^0$, and let $Z_X := Z_X'\cap Z$.

\begin{prop} \label{prop:2.1}
The projection
\[
p_2:Z\times_{{\rm Grass}(1,n)} \sU\to\P^n
\]
is surjective. Intuitively, a general point on $\P^n$ has
a line through it maximally tangent to
$X$.
\end{prop}
\begin{proof}  It suffices to consider geometric points.
Let $y\in \P^n\setminus  X$ be a geometric point. There is a linear
transformation of $\P^n$ such that the equation of $X$ is
$x_0^d +x_0^{d-1}f_1(x_1,\ldots, x_n)
+\ldots + f_d(x_1,\ldots, x_n)$, with $f_i\in  k[x_1,\ldots, x_n]$
homogeneous of degree $i$, and
such that $y$ has homogeneous coordinates $(1:0:\ldots:0)$.
Thus a line passing through $y$
has parametrization $(s:tu_1:\ldots:tu_n)\in \P^n$, for $(s:t)\in \P^1
$ and $(u_1:\ldots :u_n)\in \P^{n-1}$. The intersection
of this line with $X$ has equation
$s^d+s^{d-1}tf_1(u_1,\ldots,u_n)+\ldots + t^df_d(u_1,\ldots, u_n)$
and its intersection with $X$ will be $d$-tangent if and only
  if this equation
has the shape $(s+tu)^d$ with $(u:u_1:\ldots : u_n)\in \P^n$.
This is equivalent to the $d$ homogeneous equations
$f_i=\binom{d}{i} u^i, i=1,\ldots, d$ in $(u:u_1:\ldots u_n)\in \P^n$.
Since $d\le n$ there exists a homogeneous solution.
\end{proof}
\begin{ex} Let $(T_0: T_1: T_2)$ be homogeneous
coordinates on $\P^2$, and let $X:T_0T_1=0$, so
$d=n=2$. Clearly in this case $Z$
is simply the variety of lines through $(0:0:1)$, and $Z_X
\subset Z$ is the two points corresponding to the components of $X$.
\end{ex}
It follows from proposition \ref{prop:2.1} that
$\dim Z \ge n-1$. Let $\sZ \subset Z$ be a general
linear  section of dimension $n-1$, and write
  $\sQ:= \sZ \times_{{\rm Grass}(1,n)} \sU$.
Recall that an {\em alteration} $\pi: (\P, Y)\to (\P^n, X)$ of
$(\P^n,X)$ is a projective, generically finite morphism
$\pi:\P\to\P^n$, with $\P$ smooth,
such that $Y:=\pi^{-1}(X)$ is a strict normal crossing divisor.

\begin{lem}\label{lem:morph1} There exists a commutative diagram of schemes
\[
\xymatrix{
\P\ar[d]_\pi& \sP \ar[l]_f \ar[d]_g\ar[dr]^q\\
\P^n& \sQ\ar[r]_{\bar{p}}\ar[l]^{\bar{r}}& \sZ
}
\]
satisfying the following conditions:
\begin{enumerate}
\item $\sQ$ and $\sZ$ are as above, and $\bar{p}$ and $\bar{r}$ are the
natural maps. In
particular, $\bar{p}$ is a $\P^1$-bundle and $\bar{r}$ is surjective.
\item $\sP$ is irreducible and normal.
\item $f:\sP\to\P$ is projective, generically finite and surjective.
\item  $q:\sP\to\sZ$ is surjective.
\item There exists a normal crossings divisor
$Y \subset \P$ such that $\pi:(\P,Y) \to (\P^n,X)$
is an alteration.
\item There exists a divisor $A \subset \P$ such that $A$ meets
all the strata of $Y$ properly, and such
that $\sP\setminus f^{-1}(A) \to \P \setminus  A$ is finite.
\end{enumerate}
Given a surjection of projective
$k$-schemes $\sQ' \surj \sQ$, the map $g$ can be
taken to factor $\sP \to \sQ' \to \sQ$.
\end{lem}
\begin{proof}Since $k$ is infinite and $\sZ \subset Z$ is a general linear
space section, the map
$\bar r$ in 1) is surjective. Let $\pi$ be an alteration.
Then there will be an irreducible
component of $\P \times_{\P^n} \sQ$ dominating both $\sQ$ and $\P$.
  Taking $\sP$ to be the
normalization of this component gives 2), 3), and 4).
(To see the final assertion, one can replace
$\sQ'$ by a plane section and assume $\sQ' \to \sQ$ has finite degree.
Then substitute $\P
\times_{\P^n} \sQ'$ in the above.) Condition 5) comes from the work of 
de Jong
\cite{DJ}.

To prove 6), we use the following result  (\cite{BJAG}, Theorem 2.1.2):
\begin{thm} Let $Y\subset \P$ be a normal crossings divisor
in a smooth variety. Let $f: W \to \P$
be a finite type morphism, and assume that $W\setminus f^{-1}(Y)
\subset W$ is dense. Given $p: \P' \to \P$
the blowup of a face (stratum) of $Y$, let $Y' = p^*(Y)_{{\rm red}}$
be the reduced pullback, and let
$f': W' \to \P'$ be the strict transform of $f$. Then $(\P', Y', f')$
satisfy the same hypotheses as
$(\P, Y, f)$. There exists a composition of such blowups, $(\P_N,Y_N)  \to
\ldots \to (\P,Y)$ such
that the strict transform morphism
$f_N:W_N \to \P_N$ meets the faces of $Y_N$ properly, i.e. for
$Z \subset \P_N$ a face of codimension $r$, $f_N^{-1}(Z) \subset W_N$
has codimension $\ge r$.
\end{thm}
Replacing the alteration $\P \to \P^n$ with a composition
$\P_N \to \ldots \to \P \to \P^n$ and
changing notation, we may assume $f: \sP \to \P$ meets faces properly.
  Since $f$ has finite
degree, this amounts to saying that the fibre of $f$ over the generic
  point of any face is finite.
The existence of a divisor $A$ as in 6) is now clear.\end{proof}

\begin{lem}\label{lem:Unit} Let $p:Q\to B$ be a smooth projective
  morphism of
$k$-schemes, with geometrically connected fibers of dimension one.
  Let
$s_0,s_\infty:B\to Q$ be sections, take
$\tilde{t}\in H^0(Q,\sO_Q(s_\infty(B)-s_0(B)))$ and suppose
  that the rational function   $t$  on $Q$ determined by $\tilde{t}$ 
satisfies
\[
{\rm div}(t)=s_0(B)-s_\infty(B).
\]
  Let $\bar{B}\subset B$ be the closed subscheme of $B$ defined by the
equation $s_0=s_\infty$. Then the restriction
  $\bar{t}$ of $t$ to $p^{-1}(\bar{B})$ is a
unit, and
there is a unit $\bar{u}$ on $\bar{B}$ with $\bar{t}=p^*(\bar{u})$.
\end{lem}

\begin{proof} The hypotheses imply $p_*\sO_Q^\times =
\sO_B^\times$ and this continues to hold
after pullback. Smoothness of $p$ implies that $s_i(B)
\subset Q$ are Cartier divisors. We view
$\tilde{t}$ as an isomorphism $\tilde{t}:
\sO(s_\infty(B)) \cong \sO(s_0(B))$. By definition, the Cartier divisors
agree over $\bar B$, so there is a tautological identification
  $\tau: \sO(s_0(B))_{p^{-1}(\bar B)}
\cong \sO(s_\infty(B))_{p^{-1}(\bar B)}$. The composition
\[\tilde{t}\circ\tau\in \Gamma\Big(p^{-1}(\bar B),
  \text{Aut}\big(\sO(s_0(B))\big)\Big) =
\Gamma\Big(p^{-1}(\bar B),
  \sO^\times_{p^{-1}(\bar B)}\Big) \cong \Gamma(\bar B,\sO^\times_{\bar
B})
\]
yields the desired unit $\bar u$.
\end{proof}

Recall $\sZ \subset Z$ is a general linear section,
where $Z$ is a space of lines in $\P^n$ which
are maximally tangent to our given hypersurface $X$.
We have removed from $Z$ any possible
irreducible components consisting entirely of lines on $X$,
so the subset $Z_X \subset Z$ of lines
on $X$ is nowhere dense. We define
$\sZ_X = \sZ \cap Z_X$ and $\sZ^0 = \sZ \setminus\sZ_X$.
By generality, $\sZ^0 \subset \sZ$ is
dense. Define $\sQ^0 = \bar p^{-1}\sZ^0$ (resp. $\sP^0 =
\pi^{-1}(\sZ^0)$). The
$\P^1$-bundle
$\sQ^0 \to \sZ^0$ has a set-theoretic section $\bar{s}^0_\infty$
associating to a line $\ell$ the unique point in
$\ell\cap X$.

Consider the diagram
\ga{2.d1}{
\xymatrix{
\sP\ar[d]_\pi \ar[dr]_g \ar[r]_-{p}&
\sP\times_{\sZ} \sQ \ar[l]<3pt>_-{s_0} \ar[d]^{p_1}\\
\sZ& \sQ\ar[l]^{\bar{p}}
}
}
Here the section $s_0$ corresponds to the map $g$.
Similarly, the set-theoretic section $\bar{s}^0_\infty$
  gives rise to a set-theoretic section $s^0_\infty:\sP^0\to
   \sP^0\times_{\sZ} \sQ$.
By making a further blow-up of faces of $Y$, enlarging $A$ and
  changing notation (cf. the last part of lemma \ref{lem:morph1}),
  we may assume that the closure $\tilde{\sP}$ of
  $s^0_\infty(\sP^0)$ in  $\sP\times_{\sZ} \sQ$ is
  finite over $\P\setminus A$, hence finite over
  $\sP\setminus f^{-1}(A)$. Replacing $\sP$ with
$\tilde{\sP}$ and changing notation, we may assume
$s^0_\infty$ gives rise to another  section
$s_\infty: \sP \to \sP\times_{\sZ} \sQ$.
The picture is now
\ga{2.d2}{\xymatrix@R-20pt{Y \ar@{^{(}->}[d] && \\ \P &
\sP \ar[l]^f \ar[r]<7pt>^-{s_0}
\ar[r]<-7pt>_-{s_\infty} &
\sP\times_{\sZ} \sQ \ar[l]|-p \\ A \ar@{^{(}->}[u]
}
}

Let $\overline \sP \subset \sP$ be the closed subscheme where
$s_0 = s_\infty$. Then
\eq{2.3}{(\sP^0 \cap f^{-1}(Y))_{{\rm red}} \subset \overline \sP.
}
Indeed, we can check this down on $\P^n$, i.e. we can ignore
the alteration $\pi$. Points in
$\sP^0$ map to pairs consisting of a
line $\ell$ maximally tangent to $X$ but not lying on $X$,
together with a point $y \in \ell$.
The fibre $p^{-1}(\ell,y) = \{(\ell, y, z)\ |\ z\in \ell\}$.
The sections $s_0$ and $s_\infty$ are given respectively by $s_0(\ell,y) =
  (\ell, y,y)$ and
$s_\infty(\ell, y) = (\ell, y, \ell\cap X)$. Since
$Y = \pi^{-1}(X)$, we get the desire inclusion
\eqref{2.3} after alteration.

\begin{lem}\label{lem2.6} Possibly enlarging the divisor $A$
(preserving the hypothesis that $A$
meets faces of
$Y$ properly), there exists a rational function $t$ on $\sP\times_\sZ 
\sQ$ such that
\eq{2.4}{{(\rm div } \
t)\cap (f\circ p)^{-1}(\P \setminus A) =
(s_0 - s_\infty)\cap (f\circ p)^{-1}(\P \setminus A),
}
and such that further, $t|((f\circ p)^{-1}(Y))_{{\rm red}}
\cap (\sP^0\times_\sZ \sQ) \equiv 1$.
\end{lem}
\begin{proof} By assumption, the map $\sP\setminus
f^{-1}(A) \to \P\setminus A$ is finite. There are thus a finite
set of points of $\sP$ lying over generic points of faces of $Y$.

Let $L=p_*(\sO(s_\infty-s_0))$. As $R^ip_*(\sO(s_\infty-s_0))=0$  for 
$i>0$, $L$ is an invertible sheaf on $\sP$;   adding divisors to
$A$ meeting faces properly, we can assume that $L$ is trivial on 
$\sP\setminus f^{-1}A$. A generating section of $L$ thus gives a 
generating section $\tilde{t}$ of $\sO(s_\infty-s_0)$ over
$\sP\times_\sZ\sQ\setminus (fp)^{-1}(A)$. We let $t$ be the 
corresponding rational function on $\sP\times_\sZ\sQ$.

Clearly $t$ satisfies \eqref{2.4}. The fact that $t$ can be taken to be 
$\equiv 1$ on the indicated divisor follows from \eqref{2.3} and lemma 
\ref{lem:Unit}. Indeed, the lemma shows that the restriction of
$t$ comes from a unit on $f^{-1}Y\cap \sP^0$. Enlarging $A$,
this unit lifts to a unit on $\sP \setminus
f^{-1}A$. Normalizing $t$, we can assume this unit is $1$.
\end{proof}

\begin{proof}[Proof of Theorem \ref{thm:mot}.] We use lemma \ref{lem2.6}
to construct an effective cycle $\sD\in \sZ^n((\P\setminus A)\times 
U,1)$ with
\ga{}{
\text{Supp}(\sD)\cap (Y\setminus A)\times U\times S^1=\emptyset,\\
\partial(\sD)=N\cdot\Delta\in \sZ^n((\P\setminus A)\times U,0) \notag
}
for some integer $N\neq0$. This suffices to prove the theorem.

To construct $\sD$, we have the closed embedding
$i:\sP\times_\sZ\sQ\to \sP\times\P^n$. Let
\ga{2.8}{
\sP\times_\sZ\sQ^\#=i^{-1}((\sP\setminus f^{-1}(A))\times U),
}
and let
\ga{2.9}{
\Gamma^\#\subset \sP\times_\sZ\sQ^\#\times \P^1
}
be the closure of the graph of $t$. Let
\ga{2.10}{
\overline{\Gamma} := i_*\Gamma^\# \subset (\sP\setminus f^{-1}(A))\times 
U\times\P^1
}
be the image of $\Gamma^\#$, and let
\ga{2.11}{
\Gamma\subset  (\sP\setminus f^{-1}(A))\times U\times S^1;\quad S^1 = 
\Spec k[t_0,t_1]/(t_0 +t_1
-1)
}
be the pull-back of $\overline{\Gamma}$ via $(t_0,t_1)\mapsto (-t_0:t_1)$.
We let $\Gamma_*\subset  \sP\times \P^n\times \P^1$ be the closure of 
$\Gamma$.

By  lemma~\ref{lem2.6},  we have
\ml{2.12}{
\Gamma_*\cap (f^{-1}(Y\setminus A)\times \P^n\times \P^1\subset
(f^{-1}(Y\setminus A)\times X\times \P^1) \\ \cup(f^{-1}(Y\setminus 
A)\times \P^n\times \{1\}).
}
Thus
\ga{2.13}{
\Gamma_*\cap (f^{-1}(Y\setminus A)\times U\times S^1)=\0.
}
Also we have
\ga{2.14}{
\partial(\Gamma)=(f\times\id)^*(\Delta)\in \sZ^n((\sP\setminus 
f^{-1}(A))\times U,0).
}
Thus, setting
\ga{2.15}{
\sD:=(f\times\id)_*(\Gamma)\in  \sZ^n((\P\setminus A)\times U,1),
}
we have
\ga{2.16}{
\text{rest}(\sD)=0,\quad \partial(\sD)=(f\times\id)_*\circ 
(f\times\id)^*(\Delta)=N\cdot\Delta,
}
where $N=\deg(f)\neq0$. This completes the proof.
\end{proof}

\section{The proof of Theorem~\ref{thm:coh}}
This section is devoted to the proof of Theorem \ref{thm:coh}.

In what follows, we write simply $H^a(X,b)$ to denote either geometric 
\'etale cohomology (viz.
$X/\F_q,\ H^a(X,b):=H^a_{\text{\'et}}(X\times_{\F_q}
\overline \F_q, \Q_\ell)$ for $(\ell, {\rm char}(\F_q))=1$), or de
Rham cohomology $H^a_{DR}(X)$
for $X$ over a field of characteristic $0$ (thus $b$ plays no r\^ole here),
or rigid
cohomology with the Frobenius action
multiplied by $p^{-b}$ for $X$ over a  perfect field of
characteristic $p$.
On de Rham cohomology we denote by $F$ the Hodge filtration (\cite{De}).

\begin{lem}\label{lem3.1} Let $\P$ be smooth and le $A \subset \P$ be a 
divisor meeting faces of
$Y$ properly, where $Y \subset
\P$ is a normal crossing divisor. Let $s: H_A^a(\P,Y;0) \to 
H^a(\P,Y;0)$. Then
\begin{enumerate}
\item In the de Rham case, ${\rm Image}(s) \subset F^{1} H^a(\P,Y;0)$.
\item In the \'etale case, the eigenvalues of
Frobenius $F_a$ on ${\rm Image}(s)$ are all
divisible by $q$.
\item In the rigid case, the slopes of $F_a$ on ${\rm Image}(s)$
are all divisible by $q$.
\end{enumerate}
\end{lem}
\begin{proof}We have a diagram
\eq{3.1}{\xymatrix{H^{a-1}(Y,0) \ar[r] & H^a(\P,Y;0) \ar[r] & H^a(\P,0) \\
H^{a-1}_A(Y,0) \ar[r] \ar[u]^s & H^a_A(\P,Y;0) \ar[r]\ar[u]^s &
  H^a_A(\P,0)\ar[u]^s.
}
}
The assertion for the middle vertical arrow reduces to the comparable 
assertions for the left
and right hand vertical arrows. (In the de Rham
case, one must use the fact that ${\rm gr}_F$ is an exact
functor.) Then the spectral sequences
\ga{}{E_1^{s,t} = H^t_{A\cdot Y^{(s)}}(Y^{(s)},0) \Rightarrow 
H^{s+t}_A(Y,0) \\
E_1^{s,t} = H^t(Y^{(s)},0) \Rightarrow H^{s+t}(Y,0) \notag
}
reduce the problem to the case where the relative divisor $Y$ is smooth.
Thus it suffices to
consider the right hand vertical arrow.

Suppose for a while we work with
\'etale cohomology. We mimic Berthelot's method as in
\cite{E2}, Lemma 2.1. Let
$\ldots \subset A_j \subset A_{j-1} \ldots \subset A_0=A$
be a finite stratification by
closed subsets such that $A_{i-1}\setminus A_i$ is smooth.
The localization sequence
\ml{3.11}{\ldots \to H^b_{A_j}(\overline{\P},\Q_\ell)\to\\
H^b_{A_{j-1}}(\overline{\P},\Q_\ell)
\to H^a_{(A_{j-1}\setminus A_j)}(\overline{\P\setminus A_j},\Q_\ell)
\to \ldots
}
commutes with the Frobenius action.
  Therefore we may assume that both $A$ and $\P$ are smooth, but no
longer projective.
We consider an affine covering $\P=\cup_{i=0}^N U_i$. The spectral sequence
\ml{3.12}{E_2^{ab}=H^b\big(\ldots \to
H^b_{A}(\overline{U^{a-1}}, \Q_\ell)\to
H^b_{A}(\overline{U^{a}}, \Q_\ell)\to\\
H^b_{A}(\overline{U^{a+1}}, \Q_\ell\big) \to \ldots)
\Longrightarrow H^{a+b}_{A}(\overline{\P}, \Q_\ell)}
allows us to reduce to the case where $\P$ is smooth affine, $A\subset \P$
is smooth, where
$r={\rm codim}(A)\ge 1$. By purity, we have a functorial Gysin isomorphism
\ga{3.13}{H^{a-2r}(\overline{A}, \Q_\ell)\xrightarrow{{\rm Gysin}}
H^{a}_A(\overline{\P}, \Q_\ell(r))
  }
By functoriality, this commutes with Frobenius, and we know that the
eignevalues of Frobenius
acting on the left term are algebraic integers
(use duality with $H^*_c$ and e.g. corollary
5.5.3(iii)  in \cite{DeInt}).
But this is equivalent to saying that the Frobenius
eigenvalues on $H^{a}_A(\overline{\P}, \Q_\ell)$ are all divisible by $q^r$,
  finishing
the proof when $k=\F_q$.

The same sort of argument in the de Rham  case reduces us to showing
\ga{3.15}{{\rm Image}( H^{b}_{DR}(A)\ \to
\  H^{b+2r}_{DR}(\P)) \subset F^rH^{b+2r}_{DR}(\P))
}
when $A$ is a smooth codimension $r$ subvariety of the smooth affine $\P$.
  Blowing up $A$
in $\P$ we may assume that $A$ has codimension 1.
We take a smooth compactification $\P\subset \P'$ such that $A$ has a 
smooth
compactification $A\subset A'\subset \P$ and the divisor $ A'\cup
W, W=(\P'\setminus \P)$ is a normal crossing divisor. Then the Gysin
  map is the connecting
homomorphism of the residue sequence
\ga{3.16}{0\to \Omega^\bullet_{\P'}(\log (W))\to
\Omega^\bullet_{\P'}(\log (W+A'))
\to \Omega^{\bullet-1}_{A'}(\log (W\cap A'))\to 0.}
Since one has the exact subsequence
\ml{3.17}{0\to \Omega^{\ge (p+1)}_{\P'}(\log (W))\to
\Omega^{\ge (p+1)}_{\P'}(\log (W+A'))
\to \\
\Omega^{\ge p}_{A'}(\log (W\cap A'))\to 0,
}
one sees that $F^pH^b_{DR}(A)$ cobounds to
$F^{p+1} H^{b+2}_{DR}(\P)$. This finishes the proof when $k$ has
  charateristic 0. Finally,
when $k$ is perfect of characteristic $p$,
we use a similar argument, replacing \'etale cohomology by rigid
cohomology and then applying
\cite{E}, Lemma 2.1. This finishes the proof of the lemma.
\end{proof}
\begin{proof}[proof of theorem \ref{thm:coh}]$\P$ will be a smooth,
  projective variety of dimension
$n$,
$Y\subset
\P$ is a normal crossings divisor, and $U := \P\setminus Y$.
Consider the diagonal \eqref{8}.  Using the cycle class map from
motivic cohomology ({\it cf. }  section~\ref{sec:cycle}), we view
  the diagonal  as being a class in our theory
\ml{}{[\Delta] \in H^{2n}((\P,Y)\times U,n) \cong \bigoplus_{a+b=2n}
  H^b(U,b)  \otimes H^a(\P,Y;a)
\cong \\
\bigoplus_{a+b=2n} \text{Hom}(H^a(\P,Y;0), H^a(\P,Y;0)).
}
(The isomorphism on the right uses the existence of a
good theory of compactly supported
cohomology in our cohomology. Of course,
the homomorphism on the right is the identity.) By the
hypotheses of theorem \ref{thm:coh},
$[\Delta]$ dies in $H^{2n}((\P\setminus A,Y \setminus Y\cap A)\times U,n)$,
which means that the map $s$ below is onto:
\eq{}{H_A^a(\P,Y;0) \stackrel{s}{\surj} H^a(\P,Y;0)\to H^a(\P\setminus A,Y
\setminus Y\cap A;a).
}
The assertions of the theorem are now consequences of lemma \ref{lem3.1}
  above.
\end{proof}

\section{Cycle maps} \label{sec:cycle}  Let  $H^*_?$ be one of the
  cohomology theories used above:
\'etale cohomology with relevant twists and Galois action (\cite{DeInt}),
de Rham cohomology with Hodge filtration (\cite{De})
or rigid cohomology with Frobenius action (\cite{Be}).
  We explain how to define cycle  maps
\ga{4.10}{
\text{cl}^n:H^{2n}_M(T,\sI(A),\Z(n))\to H^{2n}_?(T, A),
}
natural with respect to smooth pull-back and projective push-forward,
  where $H^{2n}_?(T, A)$ denotes the theory on $T$ with compact supports
  relative to $A$; if $\bar{T}$ is a compactification of $T$,
  we have the usual compactly supported cohomology
$H^*_{?,c}(T):=H^*_?(\bar{T},\bar{T}\setminus T)$,
  which is canonically defined independent of the
  choice of $\bar{T}$.

For $T$ smooth over $k$ and $W$ a closed subset, we have the
  cohomology with supports $H^*_{?,W}(T)$. We have as well the
  relatively compact version $H^*_{?,W}(T,A)$ and the natural
commutative diagram
\ga{4.11}{
\xymatrix{
H^*_{?,W}(T,A)\ar[r]\ar[d]&H^*_{?}(T,A)\ar[d]\\
H^*_{?,W}(T)\ar[r]&H^*_{?}(T).}
}
If $W\cap A=\emptyset$, the map
\ga{4.3}{
H^*_{?,W}(T,A)\to H^*_{?,W}(T)}
  is an isomorphism. $H^*_?$ satisfies the {\em homotopy property}: the map
\ga{4.1}{
p^*:H^*_{?}(T,A)\to H^*_{?}(T\times\A^1,A\times\A^1)
}
is an isomorphism.

Let $R$ be the coefficient ring $H^0_?(k)$. We have the group of
  codimension $n$ cycles $\sZ^n(T)=\sZ^n(T,0)$. For $W\subset T$ a
  closed subset, we let $\sZ^n_W(T)\subset \sZ^n(T)$ be the subgroup
of cycles with support on $W$.  If  $\text{codim}_TW\ge n$,
we have the {\em purity isomorphism}
\ga{4.2}{
\cl^n_W:\sZ^n_W(T)\otimes R\to H^{2n}_{?,W}(T),}
which is natural with respect to maps $f:(T',W')\to (T,W)$,
$f^{-1}(W)\subset W'$, with $\text{codim}_{T'}W'\ge n$.
  Taking the limit over $W$ and forgetting supports defines the map
\ga{4.30}{
\cl^n:\sZ^n(T) \to H^{2n}_?(T)
}
with $\cl^n(f^*Z)=f^*\cl^n(Z)$ for $Z\in \sZ^n(T)$, and
  $f:T'\to T$ with $\text{\codim}_{T'}f^{-1}(\text{Supp}(Z))\ge n$.
  The maps $\cl^n$ are also natural with respect to projective
  push-forward and products.

Now let $A$ be a closed subset of some smooth $T$, and let $W$
be a closed subset of $T$ with $\text{codim}_TW\ge n$ and
$W\cap A=\emptyset$. Via the isomorphisms  \eqref{4.3},
  \eqref{4.2}, we have the isomorphism
\ga{4.7}{
\cl^n_W:\sZ^n_W(T)\otimes R\to H^{2n}_{?,W}(T,A)
}
Taking the limit of such $W$ and forgetting supports gives us the 
natural map
\ga{4.8}{
\cl^n_A(0):\sZ^n(T,\sI(A),0)\to H^{2n}_{?}(T,A).
}
Similarly, we have the natural map
\ga{4.9}{
\cl^n_A(1):\sZ^n(T,\sI(A),1)\to H^{2n}_{?}(T\times S^1,A\times S^1).
}
This yields the commutative diagram
\ga{4.100}{
\xymatrix{
\sZ^n(T,\sI(A),1)\ar[r]^-{\cl_A^n(1)}\ar[d]_{\delta_1^*-\delta_0^*}&
H^{2n}_{?}(T\times S^1,A\times S^1)\ar[d]^{\delta_1^*-\delta_0^*}\\
\sZ^n(T,\sI(A),1)\ar[r]_{\cl_A^n(0)}&H^{2n}_{?}(T,A).
}
}
By the homotopy property \eqref{4.1}, the right-hand vertical arrow
is zero, so $\cl^n_A(0)$ descends to the desired map
\ga{4.110}{
\cl^n_A:H^{2n}_M(T,\sI(A),\Z(n))\to H_?^{2n}(T,A).
}
The naturality of $\cl^n_A$ with respect to flat pull-back,
  projective pushforward and products follows from that of the cycle-classes
  with support.

\bibliographystyle{plain}

\renewcommand\refname{References}

\end{document}